# AN ANALOGUE OF SERRE'S CONJECTURE FOR GALOIS REPRESENTATIONS AND HECKE EIGENCLASSES IN THE MOD-$p$ COHOMOLOGY OF $GL(n, \mathbb{Z})$

AVNER ASH[1] AND WARREN SINNOTT

THE OHIO STATE UNIVERSITY
DEPARTMENT OF MATHEMATICS
231 W. 18TH AVE
COLUMBUS, OH 43210

e-mail: ash@math.ohio-state.edu,
sinnott@math.ohio-state.edu

## 1. INTRODUCTION

Let $p$ be a prime number and $\mathbb{F}$ an algebraic closure of the finite field $\mathbb{F}_p$ with $p$ elements. Let $n$ and $N$ denote positive integers, $N$ prime to $p$. We are interested in representations of the Galois group $G_\mathbb{Q} = \mathrm{G}al(\bar{\mathbb{Q}}/\mathbb{Q})$ into $GL(n, \mathbb{F})$, unramified at all finite primes not dividing $pN$. (We shall say "$\rho$ is unramified outside $pN$.") In this paper, "representation" will always mean "continuous, semisimple representation."

We choose for each prime $l$ not dividing $pN$ a Frobenius element $\mathrm{Frob}_l$ in $G_\mathbb{Q}$. We also fix a complex conjugation $\mathrm{Frob}_\infty \in G_\mathbb{Q}$. For every prime $q$ we fix a decomposition group $G_q$ with its filtration by its ramification subgroups $G_{q,i}$. We denote the whole inertia group $G_{q,0}$ by $T_q$.

Our aim is to make a conjecture about when such a representation should be attached to a cohomology class of a congruence subgroup of level $N$ of $GL(n, \mathbb{Z})$. Then we exhibit such evidence for the conjecture as we are able to.

Set $\Gamma_0(N)$ to be the subgroup of $SL(n, \mathbb{Z})$ consisting of those matrices whose first row is congruent to $(*, 0, \ldots, 0)$ modulo $N$. Let $S_N$ be the sub-semigroup of the integral matrices in $GL(n, \mathbb{Q})$ whose first row is congruent to $(*, 0, \ldots, 0)$ modulo $N$ and with determinant positive and prime to $N$.

We denote by $\mathcal{H}(N)$ the $\mathbb{F}$-algebra of double cosets $\Gamma_0(N) S_N \Gamma_0(N)$. It is commutative. This algebra acts on the cohomology and homology of $\Gamma_0(N)$ with any coefficient $\mathbb{F}S_N$ - module. When a double coset is acting on cohomology, we call it a Hecke operator. $\mathcal{H}(N)$ contains all double cosets of the form $\Gamma_0(N) D(l, k) \Gamma_0(N)$, where $D(l, k)$ is the diagonal matrix with

[1]RESEARCH PARTIALLY SUPPORTED BY NSF GRANT DMS-9531675

$k$ 1's followed by $(n-k)$ $l$'s, and $l$ is a prime not dividing N. We use the notation $T(l,k)$ for the corresponding Hecke operator.

**Definition 1.1:** Let $\mathcal{V}$ be an $\mathcal{H}(pN)$-module and suppose $v \in \mathcal{V}$ is an eigenclass for the action of $\mathcal{H}(pN)$ with $T(l,k)v = a(l,k)v$ for some $a(l,k) \in \mathbb{F}$ for all $k = 0, \ldots, n$ and all $l$ prime to $pN$. Let $\rho$ be a representation $\rho: G_\mathbb{Q} \to GL(n, \mathbb{F})$ unramified outside $pN$ such that

$$\sum_k (-1)^k l^{k(k-1)/2} a(l,k) X^k = det(I - \rho(\text{Frob}_l)X)$$

for all $l$ not dividing $pN$. Then we shall say that $\rho$ is attached to $v$.

For example, theorems of Eichler, Shimura, Deligne and Deligne-Serre imply that if $n = 2$ and $\mathcal{V}$ is the Hecke-module of classical modular forms mod $p$ of some positive weight $k$ then there exists a $\rho$ attached to any Hecke eigenform $v \in \mathcal{V}$. Such a $\rho$ is always "odd", i.e. $det\rho(\text{Frob}_\infty) = -1$.

Conversely, the weak form of Serre's conjecture [Se] states that given a mod-$p$ 2-dimensional odd Galois representation $\rho$ of the Galois group $G_\mathbb{Q}$ there exists a mod-$p$ modular Hecke eigenform $f$ with $\rho$ attached. In the strong form of Serre's conjecture, formulas for a weight ($\geq 2$), level, and nebentype character are given in terms of $\rho$ and it is conjectured that $f$ may be taken to have that weight, level, and nebentype character . For a survey of the current state of Serre's conjecture, see [Da],[Di] and [E].

The purpose of this paper is to explore an analogue of Serre's conjecture for general $n$, especially $n = 3$. We will make a conjecture in the "niveau 1" case and present some experimental evidence for it. B. Gross has made a similar conjecture in [G4] for finite arithmetic groups. The case of niveau greater than 1 is under study.

Our conjectures use cohomology classes instead of modular forms. By the Eichler-Shimura theorem, the modular form $f$ in Serre's conjecture can be replaced by a Hecke eigenclass in the cohomology of a congruence subgroup of $SL(2, \mathbb{Z})$. When $n = 2$ our conjecture will imply Serre's conjecture in the case of "niveau 1".

By way of background, we first state Conjecture 1.2, to which our Serre-type conjecture is a partial converse. Conjecture 1.2 is a slightly weaker version of Conjecture B in [A], restated here in a form suitable for this paper. As proved in [A3], Conjecture B is true when $n = 1, 2$: by Class Field Theory when $n = 1$ and the theorems of Eichler, Shimura and Deligne mentioned above when $n = 2$.

Let $V$ be a right $\mathbb{F}[GL(n, \mathbb{Z}/N)]$-module, assumed to be finite dimensional over $\mathbb{F}$. We view $V$ as an $\mathbb{F}S_N$-module by making the elements of $S_N$ act via their reductions mod $N$.

**Conjecture 1.2:** Suppose $\beta \in H^*(\Gamma_0(N), V)$ is an eigenclass for the action of the Hecke algebra $\mathcal{H}(pN)$. Then there exists a representation $\rho: G_\mathbb{Q} \to GL(n, \mathbb{F})$ unramified outside $pN$ attached to $\beta$.

We state our analogues to the weak and strong forms of Serre's conjecture below as Conjectures 2.1 and 2.2. Note that they include the case where $\rho$ is reducible. In fact, Serre's original strong conjecture could have been made in the reducible case also, with a few modifications necessitated by the fact that there is no modular form on $SL(2, \mathbb{Z})$ of level 1 and weight 2. For example if $N = 1$ and $\rho$ is a sum of the trivial character and the cyclotomic character mod $p$, the predicted weight has to be raised from 2 to $p + 1$.

In section 3 we discuss some theoretical evidence for Conjectures 2.1 and 2.2, i.e. cases where it can be proven they hold: (1) certain monomial $\rho$ when $n = p - 1$; (2) the tensor product of two odd 2-dimensional Galois representations; (3) symmetric squares of odd 2-dimensional Galois representations; and (4) 3-dimensional reducible $\rho$'s attached to cohomology classes that restrict nontrivially to the boundary of the Borel-Serre compactification $M$ of the locally symmetric space for the relevant congruence subgroup of $SL(3, \mathbb{Z})$.

We also have numerical evidence produced by computer. Because of the difficulty of producing Galois extensions of $\mathbb{Q}$ with large irreducible Galois groups, all of our numerical evidence for Conjectures 1.2, 2.1, and 2.2 as found in [AM], [AAC] and in this paper involves reducible $\rho$'s.

Our most interesting experimental evidence concerns $\rho$'s which are sums of an *even* two-dimensional representation and a one-dimensional representation. In these cases our conjecture gives a cohomological reciprocity law for the even two-dimensional representations. These even representations are not covered by Serre's conjecture. The corresponding cohomology classes in the cohomology of $GL(3, \mathbb{Z})$ must restrict to 0 on the Borel-Serre boundary of $M$.

Implicit in the definition of "attached" (Definition 1.1) is a choice of a normalization of the Satake isomorphism, called $\eta$ in [G2]. In our case $\eta$ is the one parameter subgroup

$$\eta(t) = \begin{pmatrix} t^{n-1} & 0 & \cdots & 0 & 0 \\ 0 & t^{n-2} & \cdots & 0 & 0 \\ \vdots & \vdots & \vdots & \vdots & \vdots \\ 0 & 0 & \cdots & 0 & 1 \end{pmatrix}$$

According to [G3] we should then expect that $\rho(\text{Frob}_\infty) = \pm\eta(-1)$. This is in fact verified in practice in the examples found in [AM],[AAC] and [A4]. We include this parity requirement in Conjectures 2.1 and 2.2. It is the analogue of the oddness condition in Serre's conjecture.

When $\rho$ is reducible, we have found experimentally that a finer parity requirement on $\rho(\text{Frob}_\infty)$ is needed. This requirement is formulated below as the "strict parity condition".

We should stress that almost all of the mod-$p$ cohomology classes for which we test the conjectures lift to torsion classes in characteristic 0 cohomology, and therefore they have no obvious connection with automorphic cohomology or automorphic representations.

We thank B. Gross for explaining his work to us, H. Jacquet for telling us of Ramakrishnan's work and B. Mazur for many helpful conversations. Thanks also to G. Allison and especially E. Conrad for continuing help with the computations. Thanks to J. Quer for his assistance with $\hat{A}_4$ - extensions.

## 2. Statement of the Serre-type conjectures

Let $\omega$ denote the cyclotomic character mod $p$ of $G_{\mathbb{Q}}$. For any non-negative integer $g$ let $V_g$ denote the right $GL(n, \mathbb{F})$-module of homogeneous polynomials of degree $g$ in $n$ variables with coefficients in $\mathbb{F}$.

Irreducible $GL(n, \mathbb{F}_p)$-modules are parametrized by "good" $n$-tuples of integers $(b_1, b_2, \ldots, b_n)$ where "good" means that $0 \leq b_1 - b_2, \ldots, b_{n-1} - b_n \leq p - 1$ and $0 \leq b_n \leq p - 2$. We denote the corresponding module by $F(b_1, b_2, \ldots, b_n)$. It is the restriction to $GL(n, \mathbb{F}_p)$ of the unique simple submodule of the dual Weyl module of $GL(n, \mathbb{F})$ with highest weight $(b_1, b_2, \ldots, b_n)$. The representation $F(b_1, b_2, \ldots, b_n)$ can always be embedded into $V_g$ for $g = b_1 + pb_2 + \cdots + p^{n-1}b_n$ [DW].

If $(a_1, a_2, \ldots, a_n)$ is any $n$-tuple of integers, we denote by $(a_1, a_2, \ldots, a_n)'$ a good $n$-tuple $(b_1, b_2, \ldots, b_n)$ which is congruent to $(a_1, a_2, \ldots, a_n)$ mod $p - 1$. If none of the consecutive differences $a_i - a_{i+1}$ is divisible by $p - 1$, $(a_1, a_2, \ldots, a_n)'$ is uniquely determined. In the ambiguous cases, one must interpret statements involving $(a_1, a_2, \ldots, a_n)'$ as saying that they are true for some choice of $(b_1, b_2, \ldots, b_n)$ as above. When $n = 2$ and $\rho$ is irreducible, the ambiguous cases are those in which a discrimination must be made between a Galois representation's being "peu ramifiée" and "très ramifiée" in Serre's sense [Se]. Since we have no numerical evidence for such cases when $n > 2$ we do not pursue this discrimination here. We shall see that this ambiguity also comes into play at level 1 in another way when $\rho$ is reducible, in Proposition 2.7.

**Conjecture 2.1:** Let there be given a representation (as always continuous and semisimple) $\rho : G_{\mathbb{Q}} \to GL(n, \mathbb{F})$ unramified outside $pN$. We assume that $\rho(\text{Frob}_\infty)$ has eigenvalues $\pm(1, -1, 1, -1, \ldots)$. Then there exists some integer $N'$ divisible only by the primes dividing $pN$ and an $\mathbb{F}S_{pN'}$ - module $V$ on which the $S_{pN'}$ - action factors through reduction modulo $pN'$ and an $\mathcal{H}(pN')$ - eigenclass $\beta \in H^*(\Gamma_0(N'), V)$ such that $\rho$ attached to $\beta$.

This is the weak form of our conjecture. Note that we allow $N'$ to be divisible by $p$ in this form. The necessity of the "parity" condition on $\rho(\text{Frob}_\infty)$ is strongly indicated by the experimental evidence below. In addition, from Theorem 4.3 of [AAC] and its proof one can list all possible Galois representations that can be attached to a Hecke eigenclass in the cohomology of $GL(3, \mathbb{Z})$ with coefficients in an $\mathbb{F}S_p$ module in which the action factors through reduction mod $p$, where $p = 5, 7$. They all satisfy the "parity" condition on $\rho(\text{Frob}_\infty)$. There are strong indications that the same is true for $p = 3, 11$.

We shall next give a stronger conjecture in which we predict precise analogues of the level, nebentype character, and weight in Serre's conjecture. We only venture to predict the "weight" when $p$ is odd and $\rho$ restricted to the inertia subgroup $T_p$ at $p$ is "niveau 1" in the sense of Serre [Se].

*The level:* For each prime $q$ dividing $N$ (in particular not equal to $p$) let $G$ denote the image of a chosen decomposition group at $q$ under $\rho$ and let $G_i$ be its ramification subgroups. Denote the cardinality of $G_i$ by $g_i$. Imitating the definition of the Artin conductor (as Serre does in [Se]) we set $f_q = \sum (g_i/g_0) dim M/M^{G_i}$ where $M$ is the vector space $\mathbb{F}^n$ viewed as a $G$-module under $\rho$. We then set $N(\rho) = \prod q^{f_q}$. Note that if $\rho \approx \sigma_1 \oplus \cdots \oplus \sigma_r$ with $\sigma_i$ irreducible for all $i$ then $N(\rho) = \prod N(\sigma_i)$.

*The nebentype character:* Again following [Se] we factor $\det \rho$ into $\epsilon \omega^d$ for some integer $d$, where $\epsilon$ is unramified at $p$. By Class Field Theory we can view $\epsilon$ also as an $\mathbb{F}^\times$-valued Dirichlet character $\epsilon : (\mathbb{Z}/N(\rho))^\times \to \mathbb{F}^\times$ of conductor dividing $N(\rho)$, and as such we call it $\epsilon(\rho)$.

Let $N = N(\rho)$. We have a homomorphism from $S_{pN}$ to $(\mathbb{Z}/N)^\times$ sending a matrix to its upper left-hand entry mod $N$. We use this to pull back $\epsilon$ to a character from $S_{pN}$ to $\mathbb{F}^\times$, and we denote by $\mathbb{F}_\epsilon$ the one-dimensional space on which $S_{pN}$ acts via $\epsilon$.

If $V$ is any $GL(n, \mathbb{Z}/p)$ module, we denote by $V(\epsilon)$ the $S_{pN}$ module $V \otimes \mathbb{F}_\epsilon$ where the action on the first factor is via reduction mod $p$ of the matrices in $S_{pN}$.

*The weight:* The natural generalization of the weight is an irreducible $GL(n, \mathbb{Z}/p)$-module $F(b_1, \ldots, b_n)$. Every such module can be embedded into $V_g$ for some $g$. When $n = 2$, the Eichler-Shimura theorem relates the cohomology of $V_g$ to the modular forms of weight $g + 2$.

When $\rho$ is reducible there is an additional requirement needed in our stronger conjecture which we shall call the "strict parity" condition. We consider $\rho$ as a given representation, *not* merely up to equivalence.

Suppose $\rho$ is isomorphic to the direct sum of irreducible representations $\sigma_1, \ldots, \sigma_k$ of degrees $d_1, \ldots, d_k$ respectively. We may assume there exists a standard Levi subgroup $L$ of $GL(n)$ of type $(d_1, \ldots, d_k)$ such that the image of $\rho$ lies in $L(\mathbb{F})$. "Standard" means that if $E = \{e_1, \ldots, e_n\}$ is the standard basis of $n$-space then there exists a partition of $E$ into $k$ parts of sizes $d_1, \ldots, d_k$ respectively such that $L$ is the simultaneous stabilizer of the $k$ subspaces each of which is spanned by the basis vectors in one part of the partition.

**Definition of the strict parity condition:** With $L$ and $\rho$ as above, we say that $\rho$ satisfies the strict parity condition if and only if $\rho(\text{Frob}_\infty)$ is conjugate inside $L(\mathbb{F})$ to the diagonal matrix $\pm diag(1, -1, 1, -1, \ldots)$.

**Conjecture 2.2:** Continue the hypotheses of Conjecture 2.1. Assume $p$ is odd and $\rho$ satisfies the strict parity condition. Set $N = N(\rho)$ and $\epsilon = \epsilon(\rho)$. Let $T$ denote the chosen inertia subgroup at $p$. Suppose that $\rho|T$ is conjugate inside $L(\mathbb{F})$ to a matrix of the form

$$\rho|T \sim \begin{pmatrix} \omega^{a_1} & * & \cdots & * & * \\ 0 & \omega^{a_2} & \cdots & * & * \\ \vdots & \vdots & \vdots & \vdots & \vdots \\ 0 & 0 & \cdots & 0 & \omega^{a_n} \end{pmatrix}$$

for integers $a_1, \ldots a_n$. Then we may take $N' = N$ and

$$V = F(a_1 - (n-1), a_2 - (n-2), \ldots, a_n)'(\epsilon).$$

Remarks: (1) If the degree $*$ of cohomology is the virtual cohomological dimension of $GL(n, \mathbb{Z})$ and if $p$ is prime to the orders of all torsion elements of $\Gamma_0(N)$ then we may replace $V$ by $V_g$ for any $g$ such that $V$ embeds into $V_g$. This is because in this case, $H^*(\Gamma_0(N'), V)$ embeds into $H^*(\Gamma_0(N'), V_g)$ as Hecke-module. For example, if $(a_1 - (n-1), a_2 - (n-2), \ldots, a_n)' = (b_1, \ldots, b_n)$, we can take $g = b_1 + b_2 p + \cdots + b_n p^{n-1}$.

We use this in our experimental testing when $n = 3$ and $* = 3 = vcd(GL(3, \mathbb{Z}))$. The torsion primes of $GL(3, \mathbb{Z})$ are 2 and 3, and we will always have $p > 3$ in our examples. We compute the cohomology of a congruence subgroup $\Gamma$ of $GL(3, \mathbb{Z})$ using a cellulation of a deformation retract of minimal dimension of the Borel-Serre compactification $M$ of the locally symmetric space $X/\Gamma$, where $X$ is the symmetric space of $SL(3, \mathbb{R})$. For testing the conjectures we concentrate on the degree 3 cohomology classes. This is because cohomology Hecke eigenclasses which restrict nontrivially to the Borel-Serre boundary $\partial M$ are known [AAC] to satisfy Conjecture 2.1 with reducible attached Galois representations. (It would still be of interest to verify Conjecture 2.2 for them.) On the other hand, those cohomology classes that restrict to 0 on $\partial M$ (the "interior part" of the cohomology) only occur in degrees 2 or 3. And the degree 2 and degree 3 interior parts are dual to each other by Lefschetz duality. So for testing the conjectures on the interior cohomology, it is enough to look at the degree 3 space. See [AAC] and [AT] for details, where the quotient of homology dual to the interior cohomology is called the "quasicuspidal homology". See also Proposition 2.8 below.

Since 3 is the top dimension of the retract of $M$ used to compute the cohomology, it is much easier to work in that dimension, and our programs were written only to compute $H^3$.

Also, we worked with the easiest coefficient systems to handle algorithmically. Hence our computer programs only compute cohomology with coefficients in $V_g$ for variable $g$.

In fact, models for the irreducible modules $F(b_1, \ldots, b_n)$ are not known in general. However, it is a result of Doty and Walker [DW] that the irreducible $GL(n, \mathbb{F}_p)$-module $F(b_1, b_2, \ldots, b_n)$ can be embedded in $V_g$ where $g = b_1 + b_2 p + \cdots + b_n p^{n-1}$. This is not always the minimal embedding degree, and certainly not the minimal degree in which $F(b_1, b_2, \ldots, b_n)$ occurs as a subquotient of $V_g$ in general. We don't know what the minimal $g$ should be in general.

Embedding in $V_g$ does cause problems when $g$ is too big, as happens with irreducible representations $\rho$. Conjecture 2.2 will predict some weight $F(a, b, c)$. By twisting we can always insure that $c = 0$. Reducible representations give us an extra degree of freedom that also lets us take $b$ to be small. That is why all our numerical examples below involve reducible $\rho$.

As pointed out to us by Tiep, usually there are workable modules of dimension smaller than that of $V_g$ for $g = a + pb$ in which we can embed a given $F(a, b, 0)$. For example, $F(a, b, 0)$ is always a submodule of $V_a \otimes V_b$. If we are willing to work with subquotients (although we are not guaranteed that the package of Hecke eigenvalues attached to an interior cohomology class with coefficients in a *subquotient* of $V$ will also appear in $H^3(V)$) then we can use the fact that $F(a, b, 0)$ is always a subquotient of $V_{a-b} \otimes V_b^*$ where the star denotes $\mathbb{F}$-dual. We hope to use these smaller models in some future calculations.

(2) Our assumption on the restriction of $\rho$ to inertia at $p$ corresponds to the "niveau 1" case of Serre's conjecture. We have no experimental evidence for higher niveaux, nor are we certain what the conjecture should look like in those cases.

(3) If the image of $\rho|T$ is not the whole upper triangular subgroup, for instance if it is contained in the diagonal subgroup, then there may be more than one ordering possible for the characters along the diagonal. One may then choose the ordering that leads to the smallest possible value of $g$. Other orderings would give "companion forms", as they are called when $n = 2$ [G3]. Unfortunately, in the cases we have at hand, the $g$ for the companions is too large for us to test their existence.

(4) Let $n = 3$. In the case where $\rho$ is the sum of a 1-dimensional and a 2-dimensional representation, the "strict parity condition" of the conjecture is based on computer examples discussed in Section 4 below. When $\rho$ is the sum of three 1-dimensional representations, or for larger $n$, our conjecture is based on extrapolating from the 1+2 case.

(5) If $\rho(\text{Frob}_\infty)$ has the eigenvalues $\pm(1, -1, 1, \ldots)$ the strict parity condition can always be obtained by conjugating $\rho$ and choosing $L$ appropriately. However, this will affect the order of the exponents of $\omega$ along the diagonal in $\rho|T$, and hence the "weight".

(6) Following Lemma 2.3, we will prove that Conjectures 2.1 and 2.2 are stable under twisting by powers of $\omega$ and under replacing $\rho$ with its contragredient.

To state Lemma 2.3, let $\chi : G_\mathbb{Q} \to GL(1, \mathbb{F})$ be a character of conductor dividing $N$ and let $\chi$ also denote the corresponding Dirichlet character $\chi : (\mathbb{Z}/N)^\times \to \mathbb{F}^\times$. For any $\mathbb{F} S_N$ - module $V$ as in Conjecture 1.2, let $V(\chi)$ denote the tensor product of $V$ and the one-dimensional $\mathbb{F} S_N$ - module on which $S_N$ acts via the determinant reduced modulo $N$ composed with $\chi$. We also let $V(i)$ denote the tensor product of $V$ and the one-dimensional $\mathbb{F}S_N$ module on which $S_N$ acts via the $i$-th power of the determinant mod $p$.

**Lemma 2.3:** A class $\beta \in H^*(\Gamma_0(N), V)$ is an $\mathcal{H}(N)$ - eigenclass with $\rho$ attached, if and only if $\beta \in H^*(\Gamma_0(N), V(\chi)(i))$ is an $\mathcal{H}(N)$ - eigenclass with $\rho \otimes \chi\omega^i$ attached.

**Proof:** Since $\Gamma_0(N)$ is in the kernel of the determinant, it sees $V$ and $V(\chi)(i)$ as the same module, so we can view $\beta$ also as a Hecke eigenclass in $H^*(\Gamma_0(N), V(\chi)(i))$. It is easy to check that this avatar of $\beta$ has $\rho \otimes \chi\omega^i$ attached. (The elements of $S_N$ *do* see the difference between $V$ and $V(\chi)(i)$.)

**Lemma 2.4:** $(F(c_1,\ldots,c_n)')(1)(\epsilon) = F(c_1+1,\ldots,c_n+1)'(\epsilon)$.

**Proof:** First assume $\epsilon = 1$. Let $F(c_1,\ldots,c_n)' = F(b_1,\ldots,b_n)$ where $(b_1,\ldots,b_n)$ is a good $n$-tuple. Then $(F(c_1,\ldots,c_n)')(1) = F(b_1+1,\ldots,b_n+1)$ if $b_n < p-2$ and $(F(c_1,\ldots,c_n)')(1) = F(b_1-p+2,\ldots,b_n-p+2)$ if $b_n = p-2$. In either case the parameters on the right hand side form a good $n$-tuple which is congruent mod $p-1$ to $(c_1+1,\ldots,c_n+1)$. Hence the right hand side, in either case, equals $F(c_1+1,\ldots,c_n+1)'$. To get the general result, just tensor both sides with $\mathbb{F}(\epsilon)$.

Remark: On the left hand side of Lemma 2.4 one can make any choice of $F(\ldots)'$ if there is an ambiguity, and then the proof tells what choice to make on the right hand side.

**Lemma 2.5:** Suppose $\rho : G_\mathbb{Q} \to GL(n, \mathbb{F})$ is attached to $\beta \in H^*(\Gamma_0(N), F(c_1,\ldots,c_n)'(\epsilon))$. Then $\rho \otimes \omega$ is attached to $\beta$ now viewed as a class in $H^*(\Gamma_0(N), F(c_1+1,\ldots,c_n+1)'(\epsilon))$.

**Proof:** By Lemma 2.3, $\rho \otimes \omega$ is attached to $\beta$ viewed as a class in $H^*(\Gamma_0(N), F(c_1,\ldots,c_n)'(1)(\epsilon))$. Then we are finished by Lemma 2.4.

**Proposition 2.6:** Conjectures 2.1 and 2.2 are stable under twisting by $\omega$.

**Proof:** Stability of Conjecture 2.1 under twisting follows immediately from Lemma 2.3. As for Conjecture 2.2, its stability follows from Lemma 2.5.

A small bit of evidence for the conjectures when $n = 2$ comes from the following Proposition.

**Proposition 2.7:** Conjectures 2.1 and 2.2 are true for reducible representations when $n = 2$.

**Proof:** By assumption $\rho$ is semisimple, so if it is also reducible, there are two characters $\psi$ and $\phi$ from $G_\mathbb{Q}$ to $\mathbb{F}^\times$ such that $\rho = \psi \oplus \phi$. We can write $\psi = x\omega^a$ and $\phi = y\omega^b$ where $x, y$ are unramified at $p$ and have prime-to-$p$ conductors $N$ and $M$ respectively. By Lemma 2.5, we may assume that $b = 0$, and we may also assume that $2 \leq a \leq p$. Conjecture 2.2 predicts that $\rho$ is attached to some eigenclass $\beta$ in $H^1(\Gamma_0(NM), F(a-1, 0)'(xy))$.

We view $x$ and $y$ as the reductions mod the prime $\pi$ above $p$ of characters $X$ and $Y$ into a suitable finite extension A of $\mathbb{Z}_p$. Letting $\chi$ denote the cyclotomic character for $p$, we have that $\psi$ and $\phi$ are the reductions of $X\chi^a$ and $Y$. Since $X, Y$ are unramified at $p$ they may be identified with A-valued Dirichlet characters having primitive conductors $N$ and $M$ respectively. From the parity condition on $\rho(\text{Frob}_\infty)$ we have that $XY(-1) = (-1)^{a+1}$.

Then by Proposition 1, page 3 of [W], there exists an Eisenstein series $F_{X,Y}$ of weight $a+1$, level $NM$ and nebentype $XY$, which is a Hecke eigenform with $X\chi^a \oplus Y$ as attached Galois representation. This Eisenstein series corresponds à la Eichler-Shimura (cf. Theorem 2.3 in [AS]) to a cohomology eigenclass in $H^1(\Gamma_0(NM), Sym^{a-1}(A^2)(XY))$. The reduction modulo $\pi$ of this class is the desired $\beta$. $\square$

Conjectures 2.1 and 2.2 are compatible with duality in the following sense. For any group $G$ and $\mathbb{F}[G]$-module $V$, let $V^*$ denote the dual module, i.e. $\text{Hom}(V, \mathbb{F})$ where $g$ in $G$ acts on a functional $f$ by $fg(v) = f(vg^{-1})$.

**Proposition 2.8:** Let $\rho$ be attached to a Hecke eigenclass $\alpha$ in $H^i(\Gamma_0(N), V(\epsilon))$, where $V$ is an irreducible $GL(n, \mathbb{Z}/p)$-module over $\mathbb{F}$. Let $\sigma$ be defined by $\sigma(g) = {}^t\rho(g^{-1})$. Then $\sigma$ is attached to a Hecke eigenclass $\beta$ in $H^i(\Gamma_0(N), V^*(\epsilon^{-1}))$.

**Proof:** If the invariants of $\rho$ are $N, \epsilon, V$, then the invariants of $\sigma$ are $N, \epsilon^{-1}, V^*$. This is a simple exercise, given the facts that $V^*$ is isomorphic to the contragredient of $V$, because $V$ is irreducible (see Lemma 4.6 of [AT]) and that the contragredient of $F(a_1, \ldots, a_n)'$ is $F(-a_n, \ldots, -a_1)'$.

Let $\kappa$ be the character of $\mathcal{H}$ supported by $\alpha$. Define $\lambda$ to be the character of $\mathcal{H}$ given by $\lambda(\Gamma s \Gamma) = \kappa(\Gamma s^{-1} \Gamma)$. Let $m$ be the matrix $\text{diag}(N, 1, 1)$. Consider the outer automorphism of $GL(3)$ that sends $x$ to $m^{-1} \cdot {}^t x^{-1} \cdot m$. It preserves $\Gamma_0$ and induces an isomorphism of $\mathbb{F}$-vector spaces from $H^i(\Gamma_0(N), V(\epsilon))$ to $H^i(\Gamma_0(N), V^*(\epsilon^{-1}))$. It sends $\alpha$ to a Hecke eigenclass $\beta$ supporting $\lambda$. It is a straightforward but amusing exercise to show that $\sigma$ is attached to $\beta$. $\square$

Remark: Using Lefschetz duality (as in the proof of Theorem 5.1 in [AT]) provided $\alpha$ is an interior class, one can also show the existence of a Hecke eigenclass $\gamma$ in $H^j(\Gamma_0(N), V^*(\epsilon^{-1}))$ supporting $\lambda$ where $i + j = n(n+1)/2 - 1$ is the dimension of the locally symmetric space $X/\Gamma_0(N)$.

For absolutely irreducible representations when $n = 2$ we have the following:

**Proposition 2.9:** Conjectures 2.1 and 2.2 are compatible with Serre's conjecture for irreducible representations when $n = 2$.

**Proof:** This is easily checked using the fact that $F(a, b)$ embeds into $V_{a+pb}$, and the Eichler-Shimura isomorphism between cohomology with coefficients in $V_g$ and modular forms of weight $g + 2$.

## 3. Theoretical evidence

First let $n = p - 1$. Conjecture 2.1 is true for the following type of $\rho$, as proved in [A4]. Let $K$ be the cyclotomic extension of $\mathbb{Q}$ obtained by adjoining a primitive $p$-th root of unity. Let $\theta$ denote an $\mathbb{F}$-valued ray class character on $K$ viewed as a homomorphism from $G_K$ to $\mathbb{F}^\times$. Then let $\rho : G_\mathbb{Q} \to GL(n, \mathbb{F})$ denote the induced representation from $\theta$, where $n = p - 1$.

The "weight", level and nebentype character that come in also support Conjecture 2.2. If we assume that $\rho|T$ has the stated form in Conjecture 2.2, it follows from Lemma 5.3 of [A4] that $\rho|T \equiv diag(\omega^{n-1}, \ldots, \omega, 1)$ so the predicted weight is $F(0, \ldots, 0)'$. Indeed, at the end of section 5 of [A4] it is shown that $\rho$ is attached to a cohomology class with trivial coefficients in this case.

As for the level and nebentype character appearing in [A4], they are consistent with Conjecture 2.2, but cannot be compared exactly because a different kind of congruence subgroup is used in that paper, rather than $\Gamma_0(N)$. (This is because in [A4] we work above the virtual cohomological dimension, and the congruence subgroup we use must have $p$-torsion.)

Next, let $n = 4$. Let $f_1, f_2$ be classical holomorphic modular newforms of even weights $k_1 > k_2$ for $SL(2, \mathbb{Z})$. Let $\pi_i$ be the corresponding automorphic representations on $GL(2, \mathbb{A})$, where $\mathbb{A}$ denotes the adèle group of $\mathbb{Q}$. By the work of D. Ramakrishnan [R] one knows there exists an automorphic representation $\Pi$ on $GL(4, \mathbb{A})$ which lifts the tensor product $\pi_1 \otimes \pi_2$ from $GL(2) \times GL(2)$ to $GL(4)$.

Each $f_i$ corresponds to a cohomology Hecke eigenclass in $H^1(SL(2, \mathbb{Z}), F(k_i - 2, 0))$ with a mod $p$ Galois representation $\sigma_i$ attached to it. From pp. 112-3 of [Cl] one sees that a certain Tate twist of $\Pi$ has an infinity type with $(\mathfrak{g}, K)$-cohomology with a certain finite dimensional coefficient system $E$ over $\mathbb{C}$. Taking a lattice $\Lambda \subset E$ and reducing mod $p$, we get the existence of a cohomology Hecke eigenclass $\beta$ in $H^*(SL(4, \mathbb{Z}), \Lambda/p\Lambda)$ which has attached to it the Galois representation $\rho = \sigma_1 \otimes \sigma_2 : G_\mathbb{Q} \to GL(4, \mathbb{F})$. This verifies Conjecture 2.1 for $\rho$, and the same reasoning applies when the $f_i$ have nontrivial level and nebentype.

However, in the case we are discussing we have gone further and verified Conjecture 2.2, i.e. we have determined that the weight is correctly predicted in case the $\sigma_i$ are niveau 1 (so that $\rho$ also has niveau 1). To be more precise, we have checked that the weight is correct up to a twist – to determine the twist would require following through all the normalizations involved in the parametrizations of automorphic representations by maps into the $L$-groups of $GL(2)$ and $GL(4)$. We also have to assume that $k_1$ and $k_2$ are small compared with $p$ so that $\Lambda/p\Lambda$ is irreducible. Interestingly, the requirement that $k_1$ and $k_2$ be distinct is necessitated by the fact otherwise $\Pi_\infty$ is not "regular" in the sense of [Cl] and therefore is not known to be connected with cohomology.

For the rest of this section and for the remainder of this paper we set $n = 3$. We conclude this section by discussing two series of examples that support Conjecture 2.1. The first applies to symmetric squares lifts from $GL(2)$ to $GL(3)$. The second example applies to reducible representations of type 1+2 where the 2-dimensional representation is odd. This type of representation we expect to be attached to cohomology classes that come from the boundary of the symmetric space. A fuller discussion of these examples, as of the preceding $GL(4)$ example, will be more appropriate for a separate paper.

(i) Symmetric squares: Suppose $\pi$ is a cuspidal irreducible automorphic representation for $GL(2, \mathbb{A})$ generated by some holomorphic newform $f$ of integer weight $k$ greater than 1, level $N$ and nebentype $\epsilon$. Thus by a theorem of Eichler and Shimura, $f$ corresponds to a cohomology Hecke eigenclass in $H^1(\Gamma_0(N), V_{k-2}(\epsilon))$. Let $\sigma$ denote the attached Galois representation $G_\mathbb{Q} \to GL(2, \mathbb{F})$ for a suitable $\mathbb{F}$.

Let $A : GL(2) \to GL(3)$ be the "symmetric squares" homomorphism, i.e. the adjoint representation of $GL(2)$ on the 2-by-2 matrices of trace 0 multiplied by the determinant. Let $\Pi$ be the symmetric square lifting of $\pi$ (due to Gelbart and Jacquet [GJ]). As discussed in [AS] and [AT], the automorphic representation $\Pi$ for $GL(3, \mathbb{A})$ has cohomological infinity type. That is, $\Pi_\infty$ has $(\mathfrak{g}, K)$-cohomology with a well-defined finite dimensional coefficient system $E$ over $\mathbb{C}$. Taking a lattice $\Lambda \subset E$ and reducing mod $p$, this allows one to deduce the existence of a cohomology Hecke eigenclass $\beta_f$ in $H^3(\Gamma_0(N'), (\Lambda/p\Lambda)(\epsilon'))$ for suitable $N'$ and $\epsilon'$ which has attached to it the Galois representation $G_\mathbb{Q} \to GL(3, \mathbb{F})$ given by $A \circ \sigma$.

Now do this in reverse: suppose we are given $\rho : G_\mathbb{Q} \to GL(3, \mathbb{F})$ which is the symmetric square of an irreducible 2-dimensional representation $\sigma$ and satisfies the hypotheses of Conjecture 2.1. So $\rho = A \circ \sigma$. In particular, $\sigma$ must be odd. Then by Serre's conjecture, $\sigma$ is attached to some $f$ mod $p$, and so by the preceding paragraph, $\rho$ will be attached to the corresponding $\beta_f$. Thus Conjecture 2.1 holds for $\rho$ if Serre's conjecture holds for $\sigma$.

In [AT] we show that if $\sigma$ is "niveau 1" and if $0 < k - 2 < (p-1)/2$, or if $k - 2 = (p-1)/2$ and $p = 29, 37$, or $41$, then $\Lambda/p\Lambda$ can be replaced by the weight predicted for $\rho$ by Conjecture 2.2, which is a twist of $F(2k-2, k-2, 0)$.

(ii) Borel-Serre boundary: Let $X$ be the symmetric space for $SL(3, \mathbb{R})$ and $M$ the Borel-Serre compactification of $X/\Gamma$ for an arithmetic subgroup $\Gamma$ of $SL(3, \mathbb{Z})$. Assume that $p > 3$. Then the cohomology of $\Gamma$ with trivial coefficients is canonically isomorphic to the cohomology of $M$, and one can consider the restriction map of cohomology induced by the inclusion of the boundary of $M$ into $M$.

From [LS] one sees that the cohomology of the boundary of $M$ can be given in terms of classical modular forms of weights 2 and 3 and Dirichlet characters. Therefore those boundary classes which are Hecke eigenclasses, and hence also those Hecke eigenclasses in the cohomology of $\Gamma$ that restrict nontrivially to the boundary, have attached $p$-adic Galois representations that are reducible. (Compare the statement and proof of Theorem 3.1 in [AAC] which states that any Hecke eigenclass in the cohomology of the boundary with any admissible mod $p$ coefficient module has an attached reducible Galois representation.) We can reduce the $p$-adic representation mod $p$ to get an attached $\rho$.

Remark: In the case where $\rho$ is a sum of irreducible representations of dimensions 1 and 2 respectively, the 2-dimensional representation must be odd, since it is coming from a classical modular form. Therefore, to conform with Conjecture 2.2 in predicting the weight, we must take the strict parity condition into account, although there is no parity condition on the 1-dimensional character. It remains to be seen if one can prove that every sum of a 1-dimensional and an odd 2-dimensional representation, and every sum of three 1-dimensional representations, that obeys the strict parity condition is attached to a boundary cohomology class of the predicted weight, nebentype, and level.

## 4. Computer generated examples

### Examples with $N > 1$

In [AM] we computed a number of cohomology Hecke eigenclasses with their Hecke eigenvalues for $l = 2, 3, \ldots, 97$ for the congruence subgroups $\Gamma_0(N)$ of $SL(3, \mathbb{Z})$ for prime $N$ up to about 250, trivial nebentype character and trivial coefficients $\mathbb{Z}/p$. In some of these cases we were able to find a Galois representation $\rho$ that appeared to be attached to the eigenclass in the sense that the Hecke polynomial for every $l \leq 97$ equaled the characteristic polynomial of $\mathrm{Frob}_l$.

In these cases, $\rho$ was reducible. (Cases in which the data predicted an irreducible $\rho$ had the image of $\rho$ so large that we were unable to find it.) In each of these cases we have gone back now and checked that $\rho(\mathrm{Frob}_\infty)$ has eigenvalues $(1, -1, 1)$ and that $\rho$ restricted to inertia at $p$ has eigencharacters along the diagonal $(\omega^2, \omega, 1)$ in conformity with Conjecture 2.2, in particular the strict parity condition.

These cases include the examples where $p = 3$ and the image of the 2-dimensional component is $\hat{A}_4$, and the example where $p = 5$ and the $\rho$ is a sum of three characters. In these cases, the predicted level and nebentype characters can also be checked to be correct.

One additional case that did not appear in [AM] deserves mention. There is a totally real $\hat{A}_4$ extension of $\mathbb{Q}$ unramified outside 277. This leads to a test of Conjecture 2.2 for $p = 277$ and $N = 1$ as discussed below. However, $\hat{A}_4$ is isomorphic to $SL(3, \mathbb{Z}/3)$, so the same extension gives rise to $\sigma : G_\mathbb{Q} \to SL(3, \mathbb{F})$ where now $p = 3$. Set $\rho = \sigma \oplus \omega$ with $L$ equal to the intersection of the stabilizers of the spaces of row vectors $(0, *, 0)$ and column vectors $^t(0, *, 0)$. Then the strict parity condition is satisfied if we put $\omega$ in the "middle". Then $\rho|T$ is conjugate in $L(\mathbb{F})$ to the $\mathrm{diag}(1, \omega, 1) = \mathrm{diag}(\omega^2, \omega, 1)$. The predicted weight is thus the trivial module $F(0, 0, 0)$. It is easy to see that the predicted level is $N = 277$ and the predicted nebentype character is trivial. We asked Mark McConnell to run his programs to find the Hecke module $H^3(\Gamma_0(277), \mathbb{Z}/3)$. Indeed, there was exactly one interior class (up to scalar multiples) and its Hecke eigenvalues (for primes $l \neq 3, l \leq 97$) confirmed Conjecture 2.2 in this case.

We have included an example of a calculation showing the equality of a Hecke polynomial at $l$ and the characteristic polynomial of $\text{Frob}_l$ for $p = 277, g = 90, l = 5$ at the end of the discussion of the $\hat{A}_4$ cases in the next section.

## Examples with $N = 1$

Suppose that we have an irreducible representation $\sigma : G_{\mathbb{Q}} \to GL(2, \mathbb{F})$ with the following properties:

(1) $\sigma$ is unramified outside $p$.
(2) The image of $\sigma$ has order relatively prime to $p$.
(3) $\sigma(\text{Frob}_\infty)$ is central, where $\text{Frob}_\infty$ is a complex conjugation in $G_{\mathbb{Q}}$.
(4) $\sigma(T_p)$ has order dividing $p - 1$.

Because of condition (3), Serre's Conjecture does not apply to $\sigma$. However, if we let $\rho = \omega^j \sigma \oplus \omega^k$, for suitably chosen integers $j$ and $k$, we obtain a representation to which Conjectures 2.1 and 2.2 apply, with $n = 3$ and $N = 1$: precisely, if $\sigma(\text{Frob}_\infty) = 1$, we should choose $j$ and $k$ to have opposite parity, while if $\sigma(\text{Frob}_\infty) = -1$, we should choose $j$ and $k$ to have the same parity. In this section we consider numerical evidence for Conjectures 2.1 and 2.2 for representations $\rho$ obtained in this way.

The prescription in 2.2 tells us that we should replace $\rho$ by a conjugate representation corresponding to the embedding of $GL(2) \times GL(1)$ into $GL(3)$ whose image $L$ is the Levi subgroup of the form

$$\begin{pmatrix} * & & * \\ & * & \\ * & & * \end{pmatrix};$$

then we will have

$$\rho(\text{Frob}_\infty) \sim_L \pm \begin{pmatrix} 1 & & \\ & -1 & \\ & & 1 \end{pmatrix}$$

and

$$\rho|T_p \sim_L \begin{pmatrix} \omega^a & & \\ & \omega^b & \\ & & \omega^c \end{pmatrix},$$

for certain integers $a, b, c$. In these expressions "$\sim_L$" represents conjugacy in $L$, so that $b \equiv k \bmod p-1$, while $a \bmod p-1$ and $c \bmod p-1$ can be found by examining $\omega^j \sigma | T_p$. Conjecture 2.2 then predicts that we should find a Hecke eigenclass corresponding to $\rho$ in $H^*(SL(3, \mathbb{Z}), F(a-2, b-1, c)')$; by Remark (1) after Conjecture 2.2, there should be a Hecke eigenclass corresponding to $\rho$ in $H^3(SL(3, \mathbb{Z}), V_g(\mathbb{F}_p))$, where $g = a - 2 + (b-1)p +$ $cp^2$ (assuming that $a, b, c$ have been chosen so that $(a-2, b-1, c)$ is a "good" triple). The existence of an eigenclass for the predicted value of $g$ is already a partial confirmation of the conjecture; we can usually go further and compare the characteristic polynomials of $\rho(\mathrm{Frob}_l)$ with the characteristic polynomials predicted by Conjecture 2.1, for a number of small primes $l$, and find additional confirmation in this way that Conjectures 2.1 and 2.2 are correct. The computations of cohomology and Hecke action use the programs described in [AAC].

Let $\tilde{\sigma} : G_{\mathbb{Q}} \to PGL(2, \mathbb{F})$ be the projective representation corresponding to $\sigma$, let $G$ be the image of $\tilde{\sigma}$, and let $K$ be the fixed field of the kernel of $\tilde{\sigma}$. Then $K$ will be a Galois extension of $\mathbb{Q}$, unramified outside of $p$, with Galois group isomorphic to $G$. Furthermore, $K$ is totally real, since $\tilde{\sigma}(\mathrm{Frob}_\infty) = 1$.

Since $G$ has order prime to $p$, it can be lifted to characteristic 0; hence $G$ can be realized as an irreducible subgroup of $PGL(2, \mathbb{C})$ and so is isomorphic to $A_5$, $S_4$, $A_4$, or a dihedral group.

Thus we can organize our search for representations $\sigma$ satisfying conditions (1)–(4) above by reversing these steps, as follows:

Choose a group $G$ from the list above, and search for totally real Galois extensions $K/\mathbb{Q}$ unramified outside a single prime $p \nmid \#G$ with $\mathrm{Gal}(K/\mathbb{Q}) \simeq G$. Given such a $K$, choose a projective representation $\tilde{\sigma} : G_{\mathbb{Q}} \to PGL(2, \mathbb{F})$ whose kernel has $K$ as fixed field, and consider the problem of lifting $\tilde{\sigma}$ to a representation $\sigma : G_{\mathbb{Q}} \to GL(2, \mathbb{F})$ unramified outside $p$. The following lemma (and its proof) is due to Serre; we heard about it from Richard Taylor:

**Lifting Lemma** : Let $\tilde{\sigma} : G_{\mathbb{Q}} \to PGL(2, \mathbb{F})$ be a continuous projective representation which is totally real (i.e. $\tilde{\sigma}(\mathrm{Frob}_\infty) = 1$) and unramified outside the prime $p$. Then there is a unimodular lifting $\sigma : G_{\mathbb{Q}} \to SL(2, \mathbb{F})$ of $\tilde{\sigma}$ which is unramified outside $p$.

**Proof:** We show first that a unimodular lifting exists; then we investigate the ramification. Since $\mathbb{F}$ is algebraically closed, the map $SL(2, \mathbb{F}) \to PGL(2, \mathbb{F})$ is surjective, and its kernel is $\pm 1$; hence the obstruction to lifting $\tilde{\sigma} : G_{\mathbb{Q}} \to PGL(2, \mathbb{F})$ is the class of $\tilde{\sigma}$ in $H^2(G_{\mathbb{Q}}, \pm 1)$. From class field theory we have a short exact sequence of local and global Brauer groups

$$0 \to H^2(G_{\mathbb{Q}}, \bar{\mathbb{Q}}^\times) \to \bigoplus_v H^2(G_{\mathbb{Q}_v}, \bar{\mathbb{Q}}_v^\times) \to \mathbb{Q}/\mathbb{Z} \to 0,$$

which yields, taking the kernel of multiplication by 2, the following exact sequence:

$$0 \to H^2(G_{\mathbb{Q}}, \pm 1) \to \bigoplus_v H^2(G_{\mathbb{Q}_v}, \pm 1) \to \frac{1}{2}\mathbb{Z}/\mathbb{Z}.$$

Here $v$ runs over all the places of $\mathbb{Q}$ and $G_{\mathbb{Q}_v} = \mathrm{Gal}(\bar{\mathbb{Q}}_v/\mathbb{Q}_v)$ is identified with a decomposition group of $v$ in $G_{\mathbb{Q}}$. Hence it suffices to show that the restriction $\tilde{\sigma}_v$ of $\tilde{\sigma}$ to $G_{\mathbb{Q}_v}$ can be lifted for each $v$.

If $v = \infty$, then $\tilde{\sigma}_v$ is the trivial representation, since $\tilde{\sigma}(\mathrm{Frob}_\infty) = 1$, and so $\tilde{\sigma}_v$ can be lifted to the trivial representation.

If $v$ is a finite prime not equal to $p$, then $\tilde{\sigma}_v$ factors through the maximal unramified extension of $\mathbb{Q}_v$. Therefore $\tilde{\sigma}_v$ can be lifted to an unramified representation of $G_{\mathbb{Q}_v}$ by choosing an element of $SL(2, \mathbb{F})$ which maps onto $\tilde{\sigma}_v(\mathrm{Frob}_v)$, where $\mathrm{Frob}_v$ is a Frobenius automorphism in $G_{\mathbb{Q}_v}$.

It follows that $\tilde{\sigma}_p$ must also lift, by the short exact sequence above. Hence there is a representation $\sigma : G_{\mathbb{Q}} \to SL(2, \mathbb{F})$ which lifts $\tilde{\sigma}$. It is clear that $\sigma$ is uniquely determined up to multiplication by a quadratic character of $G_{\mathbb{Q}}$.

The image of $\sigma$ is not isomorphic to $G$, since $A_5$, $S_4$, and $A_4$ have no faithful 2 dimensional representations, and the dihedral groups have no faithful unimodular 2 dimensional representations. Thus the fixed field of the kernel of $\sigma$ is a quadratic extension of $K$ and $\sigma$ may be ramified at other primes besides $p$. However, multiplying $\sigma$ by a suitable quadratic character of $G_{\mathbb{Q}}$ will produce a unimodular lifting unramified outside $p$.

This can be seen as follows. Suppose that $q \neq p$ is a prime which is ramified for $\sigma$. Then $\sigma(T_q)$ is the subgroup $\pm I$ of $SL(2, \mathbb{F})$, where $I$ is the $2 \times 2$ identity matrix. Let $G_q$ be the decomposition group of $q$ in $G_{\mathbb{Q}}$ containing $T_q$. Since $\sigma(G_q)/\sigma(T_q)$ is cyclic and $\sigma(T_q)$ is central, $\sigma(G_q)$ is abelian. Let $N$ be the commutator subgroup of $G_q$; then $N \subseteq T_q$, and $T_q/N$ may be identified with $\mathrm{Gal}(\mathbb{Q}_q(\mu_{q^\infty})/\mathbb{Q}_q)$, which in turn can be identified with $\mathrm{Gal}(\mathbb{Q}(\mu_{q^\infty})/\mathbb{Q})$. (Here $\mu_{q^\infty}$ is the group of roots of unity of order a power of $q$.) Hence there is a quadratic character $\chi_q$ of $G_{\mathbb{Q}}$, unramified outside $q$, such that
$$\sigma(\tau) = \chi_q(\tau) I, \text{ for } \tau \in T_q.$$

(A simpler argument is available for $q$ odd: then $\sigma$ is tamely ramified at $q$, and $T_q$ has a unique character of order 2, so that we may take $\chi_q$ to be the quadratic character of conductor $q$.) Let
$$\chi = \prod_q \chi_q,$$
the product taken over all rational primes $q \neq p$ which are ramified for $\sigma$. Then $\chi\sigma$ is unramified outside $p$. □

Finally, once $\sigma$ has been found, we must check that condition (4) is satisfied (it always is in the dihedral and $A_4$ cases).

Suppose then that $\sigma$ is a unimodular representation satisfying conditions (1)–(4) above. Then $\sigma|T_p = \omega^d \oplus \omega^{-d}$ for some integer $d$. Let $\rho = \omega^j \sigma \oplus \omega^k$, where the integers $j$ and $k$ satisfy the parity condition needed to guarantee that the eigenvalues of $\rho(\mathrm{Frob}_\infty)$ are either $(1, -1, 1)$ or $(-1, 1, -1)$. Then
$$\rho|T_p \simeq \omega^{j+d} \oplus \omega^{j-d} \oplus \omega^k.$$

The "strict parity condition" then requires that we take $b \equiv k \bmod p - 1$ ($a$ and $c$ may be taken to be congruent to $j + d$ and $j - d$ in either order). In the examples that follow, we choose $j$ and $k$ so as to make $g$ small: this means that we take $j = \pm d$ and $k = 1$ or $2$ (the parity of $k$ is fixed by the choice of $j$); this allows us to take $c = 0$ and $b = 1$ or $2$.

**Examples of dihedral type** We take $G$ to be a dihedral group, so that we are looking for a Galois extension $K/\mathbb{Q}$ which is totally real and unramified outside a prime $p$, and for which $\mathrm{Gal}(K/\mathbb{Q})$ is dihedral and of order prime to $p$. $p$ must be odd since $G$ has even order. Let $k$ be the maximal abelian extension of $\mathbb{Q}$ in $K$. Then since $k$ is unramifed outside $p$, $k/\mathbb{Q}$ is cyclic; it follows that $[k : \mathbb{Q}] = 2$ and that the relative degree $[K : k]$ is odd, since $G$ is dihedral. Since $k$ is also totally real, we have $p \equiv 1 \bmod 4$ and $k = \mathbb{Q}(\sqrt{p})$.

Let $T$ be an inertia group of $p$ in $G$. Since $p$ is tamely ramified in $K$, $T$ is cyclic; since $T$ surjects onto $\mathrm{Gal}(k/\mathbb{Q})$, $T \cap \mathrm{Gal}(K/k)$ must be trivial. So $T$ has order 2, and $K/k$ is unramified at $p$.

Hence the dihedral extensions we are looking for are in 1-to-1 correspondence with the set of pairs $(p, C)$, where $p$ is a prime $\equiv 1 \bmod 4$ and $C$ is a nontrivial cyclic quotient of the class group of $\mathbb{Q}(\sqrt{p})$. (It is known that the class number of $\mathbb{Q}(\sqrt{p})$ is odd and prime to $p$, so no further restrictions on $C$ are needed; see [Sl].)

We fix such a dihedral extension $K$, and let $\tilde{\sigma} : G_\mathbb{Q} \to PGL(2, \mathbb{F})$ be a projective representation corresponding to $K$. In this case $\tilde{\sigma}$ lifts to a representation $\sigma : G_\mathbb{Q} \to GL(2, \mathbb{F})$ whose image is also isomorphic to $G$, so that $\sigma$ is again totally real and unramified outside $p$. We have $\sigma|T_p \simeq \omega^{(p-1)/2} \oplus 1$, since the image of $T_p$ under $\sigma$ has order 2 and $\det \sigma = \omega^{(p-1)/2}$ (the determinant of any of the 2 dimensional irreducible representations of $G$ is the nontrivial 1 dimensional character of $G$). So we may take $\rho = \sigma \oplus \omega$, for which we will have $a = (p-1)/2$, $b = 1$, $c = 0$, so that $g = a - 2 = (p-5)/2$.

There are 6 primes $p < 1000$ which have class number $h$ greater than 1; in each of these cases $h$ itself is prime, so we obtain a unique dihedral extension unramified outside $p$:

| $p$ | 229 | 257 | 401 | 577 | 733 | 761 |
|---|---|---|---|---|---|---|
| $h$ | 3 | 3 | 5 | 7 | 3 | 3 |
| $g$ | 112 | 126 | 198 | 286 | 364 | 378 |

The dihedral group of order $2h$ has $\phi(h)/2$ two dimensional faithful representations. Hence for $p = 229, 257, 733$, and $761$, there is a unique choice of $\sigma$, while for $p = 401$ there are two choices for $\sigma$, and for $p = 577$ there are three. Using the programs described in [AAC], we find a unique interior Hecke eigenclass for $p = 229$ in weight $g = 112$, a unique interior Hecke eigenclass for $p = 257$ in weight $g = 126$, and a unique interior Hecke eigenclass for $p = 733$ in weight $g = 364$; for $p = 401$, we find a 2-dimensional interior Hecke eigenspace in weight $g = 198$, and for $p = 577$, there is a 3-dimensional interior Hecke eigenspace in weight 286. Moreover, as stated in [AAC], the Hecke action on the eigenclasses for $p = 229$ and 257 was computed for a number of small primes $l$ ($l \leq 13$ for $p = 229$; $l \leq 19$ for $p = 257$) and in each case was found to be consistent with the representation $\rho = \sigma \oplus \omega$, ($G \simeq S_3$ in these cases, so $\sigma$ is the unique irreducible 2-dimensional representation of $G$). Finally, we have also computed the Hecke action for $l \leq 7$ for the classes for $p =$401, 577, and 733, again with results consistent with Conjecture 2.2. The final pair in the table above ($p = 761, g = 378$) is just out of reach of our computer programs at present.

**Examples of $A_4$ type** We take $G$ to be $A_4$, and look for totally real Galois extensions $K/\mathbb{Q}$, unramified outside a single prime $p$, with Galois group isomorphic to $A_4$. An analysis similar to that given in the dihedral case shows that such extensions are in 1-to-1 correspondence with the set of pairs $(p, C)$, where $p$ is prime, $p \equiv 1 \bmod 3$, and $C$ is a Galois stable quotient of the ideal class group of the cubic subfield of $\mathbb{Q}(\zeta_p)$ isomorphic to the Klein four-group. These fields can be found be looking for monic quartic polynomials $f(x) \in \mathbb{Z}[x]$ with Galois group $A_4$, four real roots, and field discriminant the square of a prime $p$. A search of the Bordeaux tables [B] finds 9 such extensions in the range $p < 1000$; the values of $p$ for these extensions are $p =$163, 277, 349, 397, 547, 607, 709, 853, and 937. For each of these primes there is a unique field $K$. Note that the ramification index of $p$ in $K$ is 3.

Let $p$ be one of these primes, and let $K$ be the associated $A_4$-extension. Let $\tilde{\sigma} : G_{\mathbb{Q}} \to PGL(2, \mathbb{F})$ be a projective representation corresponding to $K$; $\tilde{\sigma}$ is unique (up to conjugacy). By the lifting lemma above, there are exactly two liftings of $\tilde{\sigma}$ to $SL(2, \mathbb{F})$ which are unramified outside $p$; if $\sigma : G_{\mathbb{Q}} \to SL(2, \mathbb{F})$ is one these liftings, $\omega^{(p-1)/2}\sigma$ is the other.

Let $\hat{K}$ be the fixed field of the kernel of $\sigma$, and let $e$ be the ramification index of $\hat{K}$ at $p$. Then $e$ is either 6 or 3 (depending on whether $\hat{K}/K$ is ramified at $p$ or not); since $\sigma(T_p)$ is abelian, we have

$$\sigma|T_p \simeq \omega^{(p-1)/e} \oplus \omega^{-(p-1)/e}.$$

It follows that replacing $\sigma$ by $\omega^{(p-1)/2}\sigma$ switches the cases $e = 6$ and $e = 3$; hence both cases do in fact occur, and we may suppose that $\sigma$ has been chosen so that $e$ is equal to 3 and $\hat{K}/K$ is unramified at $p$.

The image of $\sigma$ is isomorphic to $\hat{A}_4$, the double cover of $A_4$. $\hat{K}$ either totally real or totally complex, and to test Conjectures 2.1 and 2.2 we need to know which.

We do this as follows. Given the field $K$, described as the splitting field of a monic quartic polynomial $f(x) \in \mathbb{Z}[x]$, we obtained from Jordi Quer an element $g \in K$, written as a polynomial with integer coefficients in two of the roots $a, b$ of $f$ (since $K = \mathbb{Q}(a, b)$, for any two roots of $f$), with the property that $K(\sqrt{g})$ is an $\hat{A}_4$ extension of $\mathbb{Q}$. Quer's programs are based on the method of Crespo [C] for constructing $\widehat{A_n}$-extensions.

**Lemma** The set of $\hat{A}_4$ extensions $L$ of $\mathbb{Q}$ containing $K$ is in one-to-one correspondence with the set of square-free integers: each such $L$ can be written as $L = K(\sqrt{gm})$ for some unique square-free integer $m$.

**Proof:** From the short exact sequence
$$1 \to \{\pm 1\} \to K^\times \to (K^\times)^2 \to 1,$$
induced by squaring, we obtain the exact sequence
$$1 \to H^1(G, (K^\times)^2) \to H^2(G, \pm 1),$$
where $G = \mathrm{Gal}(K/\mathbb{Q}) \simeq A_4$; and from
$$1 \to (K^\times)^2 \to K^\times \to K^\times/(K^\times)^2 \to 1,$$
we obtain the exact sequence
$$\mathbb{Q}^\times \to (K^\times/(K^\times)^2)^G \to H^1(G, (K^\times)^2) \to 1.$$
Splicing these together gives us an exact sequence
$$\mathbb{Q}^\times \to (K^\times/(K^\times)^2)^G \to H^2(G, \pm 1);$$
by direct calculation, one sees that an element $x \in K^\times$ which representing a nontrivial class in $(K^\times/(K^\times)^2)^G$ maps to the class in $H^2(G, \pm 1)$ of the group extension
$$1 \to \{\pm 1\} \to \mathrm{Gal}(K(\sqrt{x})/\mathbb{Q}) \to G \to 1.$$
(Here we have identified $\mathrm{Gal}(K(\sqrt{x})/K)$ with $\{\pm 1\}$). $H^2(G, \pm 1)$ has order 2, and the class of $\hat{A}_4$, as an extension of $A_4$ by $\pm 1$, corresponds to the nontrivial class in $H^2(G, \pm 1)$. Hence the element $g$ provided us by Quer lies in $(K^\times/(K^\times)^2)^G$ and maps to the nontrivial class in $H^2(G, \pm 1)$; and if $x \in (K^\times/(K^\times)^2)^G$ is any element which maps to the nontrivial class in $H^2(G, \pm 1)$, then $x$ can be written in the form $gmy^2$, where $m$ is square free integer and $y \in K$. The uniqueness of $m$ results from the fact that $\mathbb{Q}^\times \cap (K^\times)^2 = (\mathbb{Q}^\times)^2$, since $K$ has no quadratic subfields. □

Thus the extension $\hat{K}$ is of the form $K(\sqrt{gm})$ for some unique square free integer $m$; of course, we are free to modify $g$ by a square as well, if convenient.

To modify $g$ so as to find a generator for $\hat{K}/K$, we proceed as follows. First of all, since $K(\sqrt{g})$ is Galois over $\mathbb{Q}$, it follows that $g^{\sigma-1}$ is a square in $K$, for every $\sigma \in G$. Let $F = \mathbb{Q}(a)$; then $[K : F] = 3$, which is odd, so that $r = \mathrm{N}_{K/F}g$ will generate the same quadratic extension of $K$ as $g$. Next we find the largest rational integer $N$ which divides $r$ in the ring of integers $O_F$ of $F$, and set $s = r/N$. Then $K(\sqrt{s})$ will be a possibly different $\hat{A}_4$-extension.

We now observe that $K(\sqrt{s})/K$ is unramified at all primes except possibly those lying above $p$ or 2. Indeed, if $q \neq p$ is an odd rational prime, then $q$ does not divide $s$ in $O_F$; hence there is some prime $\mathfrak{q}$ of $F$ dividing $q$ which does not divide $s$. Hence $F(\sqrt{s})/F$ is unramified at $\mathfrak{q}$, and so $K(\sqrt{s})/K$ is unramified at all primes of $K$ above $\mathfrak{q}$, and therefore at all primes above $q$.

To decide whether $K(\sqrt{s})/K$ is unramified at the primes above $p$, we need to examine the divisibility of $s$ by the primes above $p$ in $F$; in fact, in all the examples we treat below, $s$ turns out to be prime to $p$. This is easily detected by calculating $\mathbb{N}_{F/\mathbb{Q}}(s)$ and seeing that this integer is not divisible by $p$.

At this point we can say that $\hat{K}$ is either $K(\sqrt{s})$ or $K(\sqrt{-s})$, since the only freedom we have left is to change $s$ to $-s$. To decide which, we have to consider ramification at 2. As with odd primes, there is a prime $\mathfrak{p}$ dividing 2 in $F$ which does not divide $s$. We have the following lemma:

**Lemma** The following are equivalent:

1. $K(\sqrt{s})/K$ is unramified at all primes of $K$ above 2.
2. $F(\sqrt{s})/F$ is unramified at $\mathfrak{p}$.
3. $s \equiv u^2 \pmod{\mathfrak{p}^2}$ for some $u \in O_F$.

**Proof:** (1) $\iff$ (2) is immediate, since $K(\sqrt{s})/\mathbb{Q}$ is Galois and $K$ is unramified at 2. To show that (3) $\Rightarrow$ (2), we may change $s$ by a square and suppose that $s \equiv 1 \pmod{\mathfrak{p}^2}$; then $F(\sqrt{s})/F$ is also generated by the roots of the polynomial $x^2 + x + (1-s)/4$, hence is unramified at $\mathfrak{p}$. (Note that $(1-s)/4$ is integral at $\mathfrak{p}$ because 2 is unramified in $F$.) To show that (2) $\Rightarrow$ (3), we first observe that $(O_F/\mathfrak{p}^2)^\times \simeq \mathbb{F}(\mathfrak{p})^\times \times \mathbb{F}(\mathfrak{p})^+$ (with $\mathbb{F}(\mathfrak{p}) = O_F/\mathfrak{p}$), from which it follows that

$$s \text{ is a square in } (O_F/\mathfrak{p}^2)^\times \iff s^{q-1} \equiv 1 \pmod{\mathfrak{p}^2},$$

where $q = \#\mathbb{F}(\mathfrak{p}) = \mathbb{N}\mathfrak{p}$. Now suppose that $F(\sqrt{s})/F$ is unramified at $\mathfrak{p}$. If $\mathfrak{p}$ splits in $F(\sqrt{s})$, then $s$ is a square in $\mathbb{F}(\mathfrak{p})$ and (3) holds. On the other hand, if $\mathfrak{p}$ remains prime, then the Artin symbol of $\mathfrak{p}$ for the extension $F(\sqrt{s})/F$ is the non-trivial automorphism $\phi$ of $F(\sqrt{s})/F$; we have the congruence

$$\phi(\sqrt{s}) = -\sqrt{s} \equiv (\sqrt{s})^q \pmod{\mathfrak{p}},$$

which implies on squaring that $s \equiv s^q \pmod{\mathfrak{p}^2}$, so that (3) holds. $\square$

The method for checking whether $s$ or $-s$ satisfies condition (3) varies from case to case. Often, the simplest thing to do is to modify $s$ by a square in $K$ so that it becomes prime to 2, then considering $s$ mod 4 (this amounts to considering condition (3) for all the primes dividing 2 at once). In any case, once we have found $t = \pm s$ such that $\hat{K} = K(\sqrt{t})$, we determine whether $\hat{K}$ is totally real or totally complex by determining whether $t$ is totally positive or totally negative.

The standard part of these calculations goes as follows, in GP-PARI (1.39). We have to supply the quartic polynomial $f(x)$, the expression for $g$ as a polynomial in two of the roots ($a$ and $b$) of $f$, and the prime $p$.

```
f=????;
galois(f)
sturm(f)
factor(discf(f))
sqrt(disc(f)/discf(f))
```

```
a=mod(x,f); ff=f/(x-a);
b=[0,1,0; 0,0,1;-coeff(ff,0),-coeff(ff,1),-coeff(ff,2)];

g=????;
r=det(g); s=r/content(r);

p=????;
mod(norm(s),p)
```

The first three commands compute the Galois group of $f$, the number of real roots of $f$, and the factorization of the discriminant of the field $F = \mathbb{Q}(a)$; these are simply to confirm that we do indeed have a totally real $A_4$-extension of $\mathbb{Q}$ ramified only at a single prime. The next command calculates the index $[O_F : \mathbb{Z}[a]]$. We view the second root $b$ as the companion matrix (with entries from $F$) of the polynomial $f(x)/(x-a)$, which is an irreducible cubic over $F$, so that the norm $r = \mathbb{N}_{K/F} g$ can be computed by taking the determinant of the matrix of $g$. We obtain $s$ by removing from $r$ the gcd of its coefficients as a polynomial (of degree at most 3) in $a$. The last command checks that $s$ is prime to $p$ (provided that `mod(norm(s),p)` is nonzero, which is true for the cases examined below).

When $[O_F : \mathbb{Z}[a]] > 1$, $s$ may still be divisible by an integer dividing this index — in the cases treated below, $s$ is sometimes still divisible by 2. This can be detected by checking to see if $s/2$ is an algebraic integer:

```
content(char(s/2,x))
```

—if $s/2$ is an integer in $F$, then the content of its characteristic polynomial will be 1, in which case we replace $s$ by $s/2$. Of course we must then test the new $s$ to see whether it is still divisible by 2.

We illustrate with $p = 277$. The field $F$ is generated by a root $a$ of the polynomial $f(x) = x^4 - x^3 - 16x^2 + 3x + 1$, and the value of $g$ obtained from Quer is $9417 + 787a - 51a^2 - 133a^3 + 1359b - 627ab - 468a^2b + 114a^3b - 345b^2 - 216ab^2 + 57a^2b^2$. The PARI commands above confirm that the splitting field of $f$ is a totally real $A_4$-extension of $\mathbb{Q}$ ramified only at 277, determine that the index $[O_F : \mathbb{Z}[a]]$ equals 4, and find that $s = -150700a^3 + 193305a^2 + 2034265a + 347661$. $s$ is an element of the field $F$ such that $K(\sqrt{s})$ is an $\hat{A}_4$ extension of $\mathbb{Q}$, and such that $K(\sqrt{s})/K$ is unramified at all primes except possibly 2 and 277. In fact $K(\sqrt{s})/K$ is unramified at 277, since $\mathbb{N}_{F/\mathbb{Q}}(s) = 10483965209607696$ is not divisible by 277.

In this case, we need still to investigate whether $s$ might be divisible by 2 (since $\mathbb{Z}[a]$ has index 4 in $O_F$): to see whether $s/2$ is integral, we calculate the content of its characteristic polynomial:

```
content(char(s/2,x))
```

yields the value `1/4`—so $s$ is not divisible by 2, and thus is *not* divisible by at least one prime of $F$ above 2.

Using PARI, the following facts are easy to discover: since $\mathbb{N}(a+1) = 16$ there is a prime $\mathfrak{p}$ above 2 which divides $a+1$, so that $\mathfrak{p}$ also divides $a-1$; the number $(a-1)^2/(a+1)$ is an algebraic integer with odd norm; hence $a+1$ is in fact divisible by $\mathfrak{p}^2$. Note that $s \equiv a^2+a+1 \mod 4$; so $s \equiv 1 \mod (4, a+1)$, and therefore $s \equiv 1 \mod \mathfrak{p}^2$. So, by the lemma, $K(\sqrt{s})$ is unramified at all primes dividing 2, and therefore $\hat{K} = K(\sqrt{s})$. Finally, since the trace of $s$ is 3775974, it follows that $s$ is totally positive, and so $\hat{K}$ is totally real.

Similar calculations in PARI with the remaining primes yield the following table:

| $p$ | 163 | 277 | 349 | 397 | 547 | 607 | 709 | 853 | 937 |
|---|---|---|---|---|---|---|---|---|---|
| $\hat{K}$ | − | + | − | − | − | + | − | − | − |

where we have used + to stand for a totally real field, and − to stand for a totally complex field.

Let $\sigma : G_\mathbb{Q} \to SL(2, \mathbb{F})$ be the unimodular representation associated to $\hat{K}$. To get a 3-dimensional representation to which Conjectures 2.1 and 2.2 apply, we form $\rho = \omega^j \sigma \oplus \omega^k$ where $j$ and $k$ have different parity if $\sigma(\mathrm{Frob}_\infty) = 1$, and $j$ and $k$ have the same parity if $\sigma(\mathrm{Frob}_\infty) = -1$. Since $\hat{K}/K$ is unramified at $p$, the ramification index of $p$ in $\hat{K}/\mathbb{Q}$ is 3. Hence $\sigma|T_p \simeq \omega^{-(p-1)/3} \oplus \omega^{(p-1)/3}$.

Suppose that $\sigma(\mathrm{Frob}_\infty) = 1$: thus $\hat{K}$ is totally real, and so $p = 277$ or 607. Then we can take $j = -(p-1)/3$ and $k = 1$. For $\rho = \omega^{-(p-1)/3} \sigma \oplus \omega$, we can then take $a = (p-1)/3$, $b = 1$, and $c = 0$, giving (according to Conjecture 2.2) $g = (p-1)/3 - 2$; so $g = 90$ if $p = 277$ and $g = 200$ if $p = 607$. Using the programs discussed in [AAC] we do in fact find a unique one dimensional interior Hecke eigenspace in $H^3(SL(3,\mathbb{Z}), V_{90}(\mathbb{F}_{277}))$ and in $H^3(SL(3,\mathbb{Z}), V_{200}(\mathbb{F}_{607}))$. Furthermore, we have computed the eigenvalues of the Hecke operators $T(l,k)$ for $l = 2, 3, 5, 7, 11, 13$ and $k = 1, 2$ and found that the characteristic polynomials of Frobenius predicted by our conjecture are those arising from the representation $\rho = \omega^{-(p-1)/3} \sigma \oplus \omega$.

We can also take $j = (p-1)/3$ and $k = 1$, yielding $a = 2(p-1)/3$, $b = 1$, $c = 0$, and so (according to Conjecture 2.2) $g = 2(p-1)/3 - 2$. So $g = 182$ if $p = 277$ and $g = 402$ if $p = 607$. The first case is small enough for calculation; we find a unique one dimensional interior Hecke eigenspace in $H^3(SL(3,\mathbb{Z}), V_{182}(\mathbb{F}_{277}))$, and can verify that Conjecture 2.1 correctly predicts the characteristic polynomials of Frobenius for $l = 2, 3, 5, 7$.

We illustrate these calculations with $p = 277$, $l = 5$, and $\rho = \omega^{-92} \sigma \oplus \omega$. We consider the Hecke operators $T(5,k)$ ($k$=0, 1, 2, 3) acting on the unique one dimensional interior Hecke eigenspace in $H^3(SL(3,\mathbb{Z}), V_{90}(\mathbb{F}_{277}))$. $T(5,0)$ is the identity, and has eigenvalue 1; $T(5,3)$ arises from the scalar matrix $diag(5,5,5)$, and has eigenvalue $5^{90}$ (mod 277). According to [AAC], the eigenvalues of the Hecke operators $T(5,1)$ and $T(5,2)$ are 122 and 251

(mod 277), respectively, so that Conjecture 2.1 predicts that (in $\mathbb{F}_{277}[X]$)

$$\begin{aligned}\det(1 - \rho(\mathrm{Frob}_5)X) &= 1 - 122X + 5 \cdot 251X^2 - 5^3 \cdot 5^{90}X^3 \\ &= 1 + 155X + 147X^2 + 251X^3\end{aligned}$$

$\rho(\mathrm{Frob}_5)$ can be determined by examining the factorization of 5 in $\hat{K}$. Recall that $K$ is the splitting field of $f(X) = 1 + 3X - 16X^2 - X^3 + X^4$ and that $F$ is the field generated by a root of $f$. Since $f(X)$ factors mod 5 into the irreducible factors $X + 2$ and $X^3 + 2X^2 + 3$, 5 factors in $F$ into a prime of degree 1 and a prime of degree 3, and from this it is easy to see that $\mathrm{Frob}_5$ is a 3-cycle in $\mathrm{Gal}(K/\mathbb{Q}) \simeq A_4$. Since $\hat{K} = K(\sqrt{s})$, with $s = -150700a^3 + 193305a^2 + 2034265a + 347661 \equiv 1 \bmod 5$, the primes above 5 in $K$ split in $\hat{K}$.

Now $\hat{A}_4$ has a unique unimodular 2-dimensional representation, so we may identify $\sigma$ with this representation; we find $tr(\sigma(\tau)) = -1$ for any element of order 3 in $\hat{A}_4$, so that for $\rho = \omega^{-(p-1)/3}\sigma \oplus \omega = \omega^{-92}\sigma \oplus \omega$, we have $\det(1 - \rho(\mathrm{Frob}_5)X) = (1 - 5^{-92}(-1)X + 5^{-184}X^2)(1 - 5X) = 1 + 155x + 147x^2 + 251x^3$.

If $\sigma(\mathrm{Frob}_\infty) = -1$, we can take $j = -(p-1)/3$ and $k = 2$, giving $a = (p-1)/3$, $b = 2$, $c = 0$, and $g = (p-1)/3 - 2 + p$; so $g = 215$ if $p = 163$; $g = 463$ if $p = 349$, etc. The first case is small enough for calculation; again we find unique one dimensional interior Hecke eigenspace in $H^3(SL(3,\mathbb{Z}), V_{215}(\mathbb{F}_{163}))$, and can verify that Conjecture 2.1 correctly predicts the characteristic polynomials of Frobenius for $l = 2, 3, 5, 7$. The second case was too large for calculation of the action of the Hecke operators, but we were able to verify that there is a unique one-dimensional interior Hecke eigenspace in $H^3(SL(3,\mathbb{Z}), V_{463}(\mathbb{F}_{349}))$.

Finally, we mention several examples that led us to the "strict parity condition" of Conjecture 2.2:

- Let $p = 163$, and let $\sigma$ be the $\hat{A}_4$ representation associated to $p$ above. We have $\sigma(\mathrm{Frob}_\infty) = -1$ and $\sigma|T_p \simeq \omega^{54} \oplus \omega^{-54}$. If we take $\rho = \omega^{-52}\sigma \oplus 1$, then $\rho|T_p \simeq \omega^2 \oplus \omega^{56} \oplus 1$; if we drop the "strict parity condition," we would be able to take $a = 56$, $b = 2$, $c = 0$, which would lead to $g = 217$. But $H^3(SL(3,\mathbb{Z}), V_{217}(\mathbb{F}_{163})) = 0$.
- Let $p = 277$, and let $\sigma$ be the $\hat{A}_4$ representation associated to $p$ above. We have $\sigma(\mathrm{Frob}_\infty) = 1$ and $\sigma|T_p \simeq \omega^{92} \oplus \omega^{-92}$. If we take $\rho = \omega^{-91}\sigma \oplus 1$, then $\rho|T_p \simeq \omega \oplus \omega^{93} \oplus 1$; if we drop the "strict parity condition," we would be able to take $a = 93$, $b = 1$, $c = 0$, which would lead to $g = 91$. But $H^3(SL(3,\mathbb{Z}), V_{91}(\mathbb{F}_{277})) = 0$.

**The remaining types** The smallest examples of types $S_4$ and $A_5$ give values of $g$ which are too large to verify at present. For example, the smallest totally real $S_4$-extension ramified at only one prime has prime discriminant $p = 2777$; in this case $\sigma$ cuts out an $\widehat{S}_4$-extension of $\hat{K}/\mathbb{Q}$ ramified only at 2777. The ramification index of $p$ in $\hat{K}$ is 4, so that $\sigma|T_p \simeq \omega^{(p-1)/4} \oplus \omega^{-(p-1)/4} = \omega^{694} \oplus \omega^{-694}$. Hence if $\hat{K}$ is totally real, we may take $\rho =$

$\omega^{694}\sigma \oplus \omega$, leading to a predicted value of $g = 1386$. If $\hat{K}$ is totally complex, the best we can do is to take $\rho = \omega^{694}\sigma \oplus \omega^2$, leading to a predicted value of $g = 4163$.